\documentclass[12pt,dvips]{amsart}
\usepackage{amsfonts, amssymb, latexsym, epsfig}
\usepackage{euler, amsfonts, amssymb, latexsym, epsfig, epic}

\newtheorem{Theorem}{Theorem}[section] 

\newtheorem{Proposition}[Theorem]{Proposition} 
\newtheorem{Lemma}[Theorem]{Lemma}

\newtheorem{Corollary}[Theorem]{Corollary}
\newtheorem{Problem}[Theorem]{Problem}
\newtheorem{Definition-Lemma}[Theorem]{Definition-Lemma}
\newtheorem{Main Conjecture}[Theorem]{Main Conjecture}
\newtheorem{Conjecture}[Theorem]{Conjecture}
\theoremstyle{remark}
\newtheorem{Example}[Theorem]{Example}
\newtheorem{Remark}[Theorem]{Remark}

\theoremstyle{plain}

\newcommand{\cellsize}{12}
\newlength{\cellsz} 
\newsavebox{\cell}
\sbox{\cell}{\begin{picture}(\cellsize,\cellsize)
\put(0,0){\line(1,0){\cellsize}}
\put(0,0){\line(0,1){\cellsize}}
\put(\cellsize,0){\line(0,1){\cellsize}}
\put(0,\cellsize){\line(1,0){\cellsize}}
\end{picture}}
\newcommand\cellify[1]{\def\thearg{#1}\def\nothing{}%
\ifx\thearg\nothing
\vrule width0pt height\cellsz depth0pt\else
\hbox to 0pt{\usebox{\cell} \hss}\fi%
\vbox to \cellsz{
\vss
\hbox to \cellsz{\hss$#1$\hss}
\vss}}
\newcommand\tableau[1]{\vtop{\let\\\cr
\baselineskip -16000pt \lineskiplimit 16000pt \lineskip 0pt
\ialign{&\cellify{##}\cr#1\crcr}}}
\newcommand{\kellsize}{30}
\newlength{\kellsz} 
\newsavebox{\kell}
\sbox{\kell}{\begin{picture}(\kellsize,\kellsize)
\put(0,0){\line(1,0){\kellsize}}
\put(0,0){\line(0,1){\kellsize}}
\put(\kellsize,0){\line(0,1){\kellsize}}
\put(0,\kellsize){\line(1,0){\kellsize}}
\end{picture}}
\newcommand\kellify[1]{\def\thearg{#1}\def\nothing{}%
\ifx\thearg\nothing
\vrule width0pt height\kellsz depth0pt\else
\hbox to 0pt{\usebox{\kell} \hss}\fi%
\vbox to \kellsz{
\vss
\hbox to \kellsz{\hss$#1$\hss}
\vss}}
\newcommand\ktableau[1]{\vtop{\let\\\cr
\baselineskip -16000pt \lineskiplimit 16000pt \lineskip 0pt
\ialign{&\kellify{##}\cr#1\crcr}}}

\begin{document}
\pagestyle{plain}

\mbox{}
\title{A jeu de taquin theory for increasing tableaux, with
applications to $K$-theoretic Schubert calculus} 
\author{Hugh Thomas}
\address{Department of Mathematics and Statistics, University of New Brunswick, Fredericton, New Brunswick, E3B 5A3, Canada }
\email{hugh@math.unb.ca}

\author{Alexander Yong}
\address{Department of Mathematics, University of Minnesota, Minneapolis, MN 55455, USA}

\email{ayong@math.umn.edu}

\date{September 16, 2007}
\thanks{HT was supported by an NSERC Discovery grant. AY 
was supported by NSF grant 0601010.}
\maketitle
\begin{abstract}
We introduce a theory of 
\emph{jeu de taquin} for \emph{increasing tableaux}, 
extending fundamental work of 
[Sch\"{u}tzenberger '77] for standard Young tableaux. We apply this to
give a new combinatorial rule for the $K$-theory Schubert
calculus of Grassmannians via \emph{$K$-theoretic jeu de taquin}, 
providing an alternative 
to the rules of [Buch '02] and others. This rule
naturally generalizes to give a conjectural root-system uniform rule 
for any minuscule flag variety $G/P$, extending [Thomas-Yong '06]. 
We also present analogues of results of Fomin, Haiman, Schensted
and Sch\"{u}tzenberger.
\end{abstract}

\tableofcontents

\section{Introduction}

In this paper, we introduce a 
\emph{jeu de taquin} type
theory for \emph{increasing tableaux}, extending
Sch\"{u}tzenberger's fundamental framework \cite{Schutzenberger}
to the ($K$-theoretic) \emph{Grothendieck polynomial}
context introduced by Lascoux and Sch\"{u}tzenberger~\cite{LS:Hopf}.

One motivation and application for this work comes from Schubert calculus.
Let $X=Gr(k,{\mathbb C}^n)$ be the Grassmannian of $k$-planes
in ${\mathbb C}^n$ and let $K(X)$ be the {\bf Grothendieck ring}
of algebraic vector bundles over $X$, see, e.g., the expositions \cite{brion:lectures, Buch:expository} 
for definitions and discussion. To each partition, as identified with its
{\bf Young shape} $\lambda\subseteq \Lambda:=k\times (n-k)$,
let $X_{\lambda}$ be the
associated Schubert variety and ${\mathcal O}_{X_{\lambda}}$ its
structure sheaf. The classes $\{[\mathcal O_{X_{\lambda}}]\}\subseteq
K(X)$ form an additive ${\mathbb Z}$-basis of $K(X)$. 
The {\bf ($K$-theoretic) Schubert structure constants} $C_{\lambda,\mu}^{\nu}$
are defined by
\[[{\mathcal O}_{X_\lambda}]\cdot [{\mathcal O}_{X_\mu}]=\sum_{\nu\subseteq \Lambda}
C_{\lambda,\mu}^{\nu}[{\mathcal O}_{X_\nu}].\]
Buch's rule~\cite{Buch:KLR} established alternation of
sign, 
i.e., $(-1)^{|\nu|-|\lambda|-|\mu|}C_{\lambda,\mu}^{\nu}\in {\mathbb N}.$ 

There has been significant interest in the Grothendieck
ring of $X$ and of related varieties, see work on, e.g.,
quiver loci \cite{Buch:original, Buch:JAMS, Miller:Duke, BKSTY}, 
Hilbert series of determinantal ideals \cite{KM, KMY, KMYII},
applications to invariants of matroids \cite{Speyer}, 
and in relation to representation
theory \cite{Griffeth.Ram, Lenart.Postnikov, Willems}. See also
work of \cite{Lam.P} concerning combinatorial Hopf algebras.

 We aim to provide 
unifying foundational combinatorics in support of further such developments.
Evidence of the efficacy of this approach is provided through our study of
minuscule Schubert calculus; other uses are also suggested. 
In particular, as a \emph{non}-algebraic geometric application, 
in forthcoming work \cite{Thomas.Yong:VI}, we relate the ideas in this paper
to \cite{BKSTY} and 
the study of longest strictly increasing 
subsequences in random words.

The classical setting for the
{\bf Littlewood-Richardson coefficients} is the {\bf cohomology case} 
when 
\[|\lambda|+|\mu|=|\nu|, \mbox{ where 
$|\lambda|=\sum_{i}\lambda_i$ is the {\bf size} of $\lambda$}.\] 
Here, 
$C_{\lambda,\mu}^{\nu}$ counts 
points in the intersection of three general 
Schubert varieties. These numbers determine the ring structure  
of the cohomology $H^{\star}(X,{\mathbb Q})$. 
Combinatorially, they are governed by the tableau theory
of Schur polynomials. Sch\"{u}tzenberger's \emph{jeu de taquin} theory,
\cite{Schutzenberger} by which the first modern statement and proof
of a Littlewood-Richardson rule was constructed, has had a central impact here.

While $H^{\star}(X,{\mathbb Q})$ reflects important geometric data about
$X$, this is even more true of $K(X)$. The combinatorics of the latter is
encoded by the Grothendieck
polynomials of Lascoux and Sch\"{u}tzenberger \cite{LS:Hopf} (for more details,
see the Appendix). 
This richer environment parallels the Schur polynomial setting,
as demonstrated by, e.g.,~\cite{Lenart, Buch:KLR, BKSTY}. However, basic
gaps in this comparison remain. In particular, one lacks an analogue of the
\emph{jeu de taquin} theory. 
This also raises questions of intrinsic combinatorial
interest.

We introduce a jeu de taquin construction, thereby allowing for
$K$-theoretic generalizations of a number of 
results from algebraic combinatorics. In particular, we give an analogue
of Sch\"{u}tzenberger's Littlewood-Richardson rule. In addition, 
we also extend
Fomin's \emph{growth diagrams}, allowing for, e.g., a generalization of
Sch\"{u}tzenberger's \emph{evacuation involution}. On the other hand, 
it is interesting that natural
generalizations of some results from the classical theory are \emph{not} true,
underlining some basic combinatorial obstructions. 

One feature
of our rule is that it has a natural conjectural 
generalization to any minuscule flag variety $G/P$,
extending our earlier work~\cite{Thomas.Yong, Thomas.Yong:DE}; this provides
the first generalized Littlewood-Richardson formula (even conjectural)
for $K$-theory, outside of the Grassmannians. 
(There are already a number of more specialized   
$K$-theoretic Schubert calculus formulas proven for any $G/P$,
such as the Pieri-type formulas of \cite{Lenart.Postnikov} and others). 

\subsection{Main definitions}
An {\bf increasing tableau} $T$ of shape $\nu/\lambda$ is a filling of 
the skew shape ${\tt shape}(T)=\nu/\lambda$
with $\{1,2,\ldots, q\}$ where $q\leq |\nu/\lambda|$ such that 
the entries of $T$ strictly increase along each row
and column. We write $\max T$ for the maximum entry in $T$.  
In particular, when $\max T=|\nu/\lambda|$ and each label appears
exactly once, $T$ is a {\bf standard Young
tableau}. Let ${\tt INC}(\nu/\lambda)$ be the set of these increasing 
tableaux 
and ${\tt SYT}(\nu/\lambda)$ be the set of standard Young tableaux
for $\nu/\lambda$. Below we give an example of 
an increasing tableau and a standard Young tableau, each
of shape $\nu/\lambda=(5,3,1)/(2,1)$:
\[\tableau{{ \ }&{ \ }&{1}&{2}&{3}\\{ \ }&{1}&{3}\\{2}}\in 
{\tt INC}((5,3,1)/(2,1)) \mbox{ \ \ \ \ \  \ \ \ \ \ } 
\tableau{{\ }&{\ }&{1}&{4}&{6}\\{ \ }&{2}&{5}\\{3}}\in {\tt SYT}((5,3,1)/(2,1))
\]
We also need to define the {\bf superstandard Young tableau} $S_{\mu}$ of shape
$\lambda$ to be the standard Young tableau that fills the first row with
$1,2,\ldots,\lambda_1$, the second row with $\lambda_1+1,\lambda_1+2,\ldots,
\lambda_1+\lambda_2$ etc. For example,
\[S_{(5,3,3,1)}=\tableau{{1}&{2}&{3}&{4}&{5}\\{6}&{7}&{8}\\{9}&{10}&{11}
\\{12}}.\]

A {\bf short ribbon} $R$ is a connected
skew shape that does not contain a $2\times 2$ subshape and where
each row and column contains at most two boxes. A 
{\bf alternating ribbon} is a filling of a short ribbon
$R$ with two symbols where
adjacent boxes are filled differently. We define 
${\tt switch}(R)$ to be the alternating ribbon of the same shape as
$R$ but where each box is instead filled with the other symbol. 
For example, we have:
\[R=\tableau{&&{\circ}&{\bullet}\\&{\circ}&{\bullet}\\{\circ}&{\bullet}\\{\bullet}}
\mbox{ \ \ \ \ \ \ \ \ \ \ \ \ \ \ 
${\tt switch}(R)=
\tableau{&&{\bullet}&{\circ}\\&{\bullet}&{\circ}\\{\bullet}&{\circ}\\{\circ}}$.} 
\]
By definition, if $R$ is a ribbon consisting of a single box, $\tt switch$ does nothing to 
it. We define $\tt switch$ to act on a skew shape consisting of multiple
connected components, each of which is a short ribbon, by acting on each 
separately.  

Our starting point is the following new idea. Given $T\in {\tt INC}(\nu/\lambda)$,
an {\bf inner corner} is any maximally southeast box $x\in\lambda$. Now fix a set
$\{x_1,\ldots,x_s\}$ of inner corners and let each of these
boxes is filled with a ``$\bullet$''.  Consider the union of short ribbons
$R_1$ which is made of boxes with entries $\bullet$ or $1$.  Apply
$\tt switch$ to $R_1$.  Now let $R_2$ be the union of short ribbons 
consisting of boxes with entries $\bullet$ or $2$, and proceed as before.
Repeat this process $\max T$ times, in other words, until the $\bullet$'s
have been switched past all the entries of $T$. 
The final placement of the numerical 
entries gives  $K{\tt jdt}_{\{x_i\}}(T)$.

\begin{Example}
Let $T=\tableau{{ \ }&{\ }&{1}&{2}&{3}\\{ \ }&{2}&{3}\\{2}}$
as above and $\{x_i\}$ as indicated below:
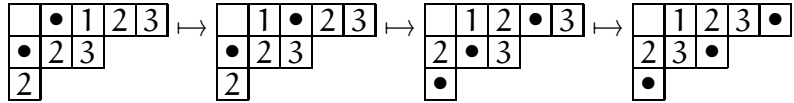
\[\tableau{{ \ }&{ \bullet }&{1}&{2}&{3}\\{ \bullet }&{2}&{3}\\{2}}\mapsto
\tableau{{ \ }&{ 1 }&{\bullet}&{2}&{3}\\{ \bullet }&{2}&{3}\\{2}}
\mapsto 
\tableau{{ \ }&{ 1 }&{2}&{\bullet}&{3}\\{ 2 }&{\bullet}&{3}\\{\bullet}}
\mapsto 
\tableau{{ \ }&{ 1 }&{2}&{3}&{\bullet}\\{ 2 }&{3}&{\bullet}\\{\bullet}}
\]
and therefore $K{\tt jdt}_{\{x_i\}}=\tableau{{ \ }&{ 1 }&{2}&{3}\\{ 2 }&{3}}$.
\end{Example}

It is easy to see that $K{\tt jdt}_{\{x_i\}}(T)$ is an increasing tableau
also. Moreover, if $T$ is a standard Young tableau, and only one corner $x$
is selected, the result is
an {\bf ordinary jeu de taquin slide} ${\tt jdt}_{x}(T)$.
Given $T\in {\tt INC}(\nu/\lambda)$ we can iterate applying $K{\tt jdt}$-slides until no such moves are possible. The result $K{\tt rect}(T)$, 
which we call
\emph{a} $K$-{\bf rectification} of $T$, is an increasing 
tableau of straight shape, i.e., one whose shape is given by some partition
$\lambda$. We will refer to the choice of intermediate $K{\tt jdt}$ slides
as a {\bf rectification order}.

\begin{Theorem}
\label{thm:welldefined}
Let $T\in {\tt INC}(\nu/\lambda)$. If $K{\tt rect}(T)$ is
a superstandard tableau $S_{\mu}$ for some rectification order, then 
$K{\tt rect}(T)=S_{\mu}$ for any rectification order.
\end{Theorem}

It will also be convenient
to define {\bf reverse slides} $K{\tt revjdt}_{\{x_i\}}(T)$ of 
$T\in {\tt INC}(\nu/\lambda)$, 
where now each $x_i$ is an {\bf outer corner}, i.e., 
a maximally northwest box $x\in \Lambda\setminus\nu$. We can
similarly define {\bf reverse rectification} $K{\tt revrect}(T)$.
Clearly, Theorem~\ref{thm:welldefined} also implies the ``reverse version''.
When we refer to {\bf slides}, we mean either $K{\tt jdt}$ or 
$K{\tt revjdt}$ operations.

This result may be compared to what is often called 
the ``confluence theorem'' or the ``First Fundamental Theorem'' in the
the original setting of \cite{Schutzenberger}. There, the superstandard
assumption is unnecessary and so rectification is always well-defined. However
this is not true in our more general context, and thus one has the task of 
recognizing the additional hypothesis needed.

\begin{Example}
\label{exa:twodiff}
Consider the following two $K$-rectifications of the same skew tableau $T$:
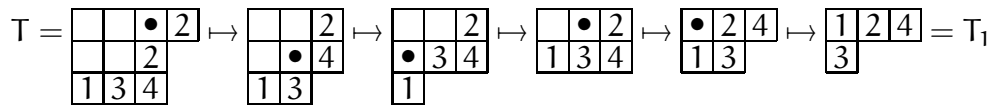
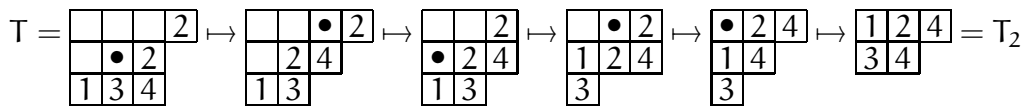
\[T=\tableau{{\ }&{\ }&{\bullet }&{2}\\{\ }&{\ }&{2}\\{1}&{3}&{4}}\mapsto 
\tableau{{\ }&{\ }&{2 }\\{\ }&{\bullet }&{4}\\{1}&{3}}\mapsto
\tableau{{\ }&{\ }&{2 }\\{\bullet }&{3 }&{4}\\{1}}\mapsto
\tableau{{\ }&{\bullet }&{2 }\\{1 }&{3 }&{4}}\mapsto
\tableau{{\bullet }&{2 }&{4 }\\{1 }&{3 }}\mapsto
\tableau{{1 }&{2 }&{4 }\\{3 }}=T_1
\]
and
\[T=\tableau{{\ }&{\ }&{\ }&{2}\\{\ }&{\bullet }&{2}\\{1}&{3}&{4}}\mapsto 
\tableau{{\ }&{\ }&{\bullet }&{2}\\{\ }&{2 }&{4}\\{1}&{3}}\mapsto 
\tableau{{\ }&{\ }&{2 }\\{\bullet }&{2 }&{4}\\{1}&{3}}\mapsto 
\tableau{{\ }&{\bullet }&{2 }\\{1 }&{2 }&{4}\\{3}}\mapsto
 \tableau{{\bullet }&{2 }&{4 }\\{1 }&{4 }\\{3}}\mapsto
 \tableau{{1 }&{2 }&{4 }\\{3 }&{4 }}=T_2
\]
Now $T_1\neq T_2$. However, neither rectification is superstandard.
\end{Example}

We need Theorem~\ref{thm:welldefined} to state our 
new combinatorial rule for $C_{\lambda,\mu}^{\nu}$:

\begin{Theorem}
\label{thm:KLR}
$(-1)^{|\nu|-|\lambda|-|\mu|}C_{\lambda,\mu}^{\nu}$ counts the
number of $T\in {\tt INC}(\nu/\lambda)$ where 
$K{\tt rect}(T)=S_{\mu}$. 
\end{Theorem}

\begin{Example}
\label{exa:witness}
The computation 
$C_{(2,2),(2,1)}^{(3,2,2,1)}=-2$ is witnessed by the increasing tableaux:
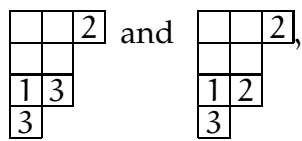
\[\tableau{{\ }&{ \ }&{2}\\{\ }&{\ }\\{1}&{3}\\{3}} \mbox{\ \ and  \ \ }
\tableau{{\ }&{ \ }&{2}\\{\ }&{\ }\\{1}&{2}\\{3}},
\]
which both rectify to $\tableau{{1}&{2}\\{3}}$. 
\end{Example}

One can replace the superstandard assumption by some other classes $\{C_\mu\}$
of tableau (most obviously the one where we consecutively number columns rather than rows), but we focus on the superstandard choice in this paper. 

We will give a self-contained proof of Theorem~\ref{thm:KLR}, once granted
Lenart's Pieri rule~\cite{Lenart}. 

A short review about past work on $K$-theoretic Littlewood-Richardson rules 
is in order:
The first rule for $C_{\lambda,\mu}^{\nu}$ 
was given by Buch~\cite{Buch:KLR}, who gave a generalization of
the \emph{reverse lattice word} formulation of the classical Littlewood-Richardson 
rule. That formula utilized the new idea of \emph{set-valued tableaux} (see the Appendix). Afterwards, another formula was given by Lascoux~\cite{lascoux:transition} in terms of counting paths in a certain tree (generalizing the 
\emph{Lascoux-Sch\"{u}tzenberger tree}, see, e.g., \cite{Manivel}). 
In \cite{KY}, Lascoux's rule was reformulated
in terms of \emph{diagram marching moves}, and it was also extended
to compute a wider class of $K$-theoretic Schubert structure constants. 
More recently, in \cite{BKSTY}, a rule was given for another class
of combinatorial numbers generalizing $C_{\lambda,\mu}^{\nu}$.
This rule specializes to a new formula for $C_{\lambda,\mu}^{\nu}$ and
in fact gives an independent proof of Buch's rule.

\subsection{Minuscule Schubert calculus}
In earlier work~\cite{Thomas.Yong, Thomas.Yong:DE}, we introduced 
root-system uniform combinatorial rules for minuscule Schubert calculus. 
Theorem~\ref{thm:KLR} has the advantage that it admits
a straightforward conjectural generalization to the minuscule setting.
We state one form of our conjecture below; more details will appear in forthcoming work. 

Let $G$ be a complex, connected reductive Lie group with root system
$\Phi$, positive roots $\Phi^{+}$ and base of simple roots $\Delta$.
To each subset of $\Delta$ is associated a parabolic subgroup $P$. The
{\bf generalized flag variety} $G/P$ has {\bf Schubert varieties}
\[X_{w}:={\overline{B_{-}wP/P}} \mbox{ \ for \ } wW_P\in W/W_P,\] 
where $W$ is the
Weyl group of $G$ and $W_P$ is the parabolic subgroup of $W$ corresponding
to $P$. Let $K(G/P)$ be the Grothendieck ring of $G/P$, with a 
basis of Schubert structure sheaves $\{[{\mathcal O}_{X_w}]\}$. Define
Schubert structure constants $C_{u,v}^{w}(G/P)$ as before, by
\[[{\mathcal O}_{X_u}]\cdot [{\mathcal O}_{X_{v}}]=\sum_{wW_P\in W/W_P}
C_{u,v}^{w}(G/P)[{\mathcal O}_{X_{w}}].\]
Brion~\cite{brion:lectures} has established that 
\[(-1)^{\ell(w)-\ell(u)-\ell(v)}C_{u,v}^{w}(G/P)\in {\mathbb N},\] 
where $\ell(w)$ is the \emph{Coxeter length} of the minimal length coset
representative of $wW_P$.

A maximal parabolic subgroup $P$ is said to be {\bf minuscule} if the 
associated fundamental weight $\omega_{P}$
satisfies $\langle\omega_{P},\alpha^{\vee}\rangle\leq 1$ for all 
$\alpha\in\Phi^{+}$ under the usual pairing between weights and coroots.
The {\bf minuscule flag varieties} $G/P$ are classified into five
infinite families and two exceptional cases (the type $A_{n-1}$ cases
are the Grassmannians $Gr(k,{\mathbb C}^n)$).

Associated to each minuscule $G/P$ is a planar poset 
$(\Lambda_{G/P},\prec)$, obtained as a subposet of the poset of
positive roots $\Omega_{G^{\vee}}$ for the dual root system of $G$. 
In this context,
{\bf shapes} $\lambda$ are lower order ideals in this poset. These shapes
are in bijection with the cosets $wW_P$ indexing the Schubert varieties; in
particular, 
if $wW_P\leftrightarrow \nu$ under this bijection, 
$\ell(w)=|\nu|$. Define a 
{\bf skew shape}
$\nu/\lambda:=\nu\setminus\lambda$ 
to be a set theoretic difference of two shapes. Define an {\bf increasing tableau}
of shape $\nu/\lambda$ 
to be an assignment 
\[{\tt label}:\nu/\lambda\rightarrow\{1,2,\ldots,q\}\] 
such that ${\tt label}(x)<{\tt label}(y)$ whenever
$x\prec y$, and where each label appears at least once. 
An {\bf inner corner} of $\nu/\lambda$ 
is a maximal element $x\in \Lambda_{G/P}$
that is below some element in $\nu/\lambda$. With these definitions, we define
notions of ${\tt INC}_{G/P}(\nu/\lambda),
K{\tt jdt}_{G/P;\{x_i\}}, K{\tt rect}_{G/P}$, superstandard $S_{\mu}$, etc., 
in a manner analogous to those we have given for the Grassmannian. The following
rule is new for all minuscule $G/P$:

\begin{Conjecture}
\label{conjecture:mainone}
For any minuscule $G/P$, $(-1)^{|\nu|-|\lambda|-|\mu|}
C_{\lambda,\mu}^{\nu}(G/P)$ equals
the the number of $T\in {\tt INC}_{G/P}(\nu/\lambda)$ such that 
$K{\tt rect}_{G/P}(T)=S_{\mu}$. 
\end{Conjecture}

Implicit in this conjecture is the conjecture that an analogue of 
Theorem~\ref{thm:welldefined} holds.
A weaker form of these conjectures
is that there is a tableau $C_{\mu}$ for each shape $\mu$ 
such that the aforementioned conjectures hold after replacing $S_{\mu}$ by~$C_{\mu}$. 

Briefly, using the ideas contained in this paper, together with those in
\cite{Thomas.Yong, Thomas.Yong:DE} it is not hard to show that
$K{\tt jdt}_{G/P;\{x_i\}}$ is well-defined. The next aim is 
to establish the analogue of 
Theorem~\ref{thm:welldefined}. Once this is achieved we can prove that
our conjectural rule defines an associative, commutative ring with an 
additive ${\mathbb Z}$-basis indexed by shapes. 
It would then remain to show that
such rules compute the correct geometric numbers. 

\subsection{Organization of this paper}
In Section~2, we introduce an analogue of Fomin's
\emph{growth diagrams}, which compute $K$-rectifications; 
their symmetries
make it possible to give a simple proof of the \emph{infusion involution}
of Section 3. In Section~4, we again exploit growth diagrams to 
give an analogue of Sch\"{u}tzenberger's \emph{evacuation involution}.
In Section~5, we use the infusion involution to show that if 
Theorem~\ref{thm:welldefined} holds, then
Theorem~\ref{thm:KLR} indeed 
computes Schubert calculus. Theorem~\ref{thm:welldefined} itself is 
actually proved in Section~6, 
where we also need a connection to 
\emph{longest strictly increasing subsequences} of reading words of tableaux.
In Section~7, we give more details of our conjectural minuscule Schubert calculus rule, together with an example.
In Section~8, we give counterexamples to natural
analogues of various results that are true for classical
Young tableau theory. Finally, in Section~9 we give some concluding remarks and further conjectures.
In order to be self-contained, we provide an appendix giving 
background about Grothendieck polynomials so that our results 
can be given a completely elementary and concrete origin.

\section{Growth diagrams}
A construction that is important to this paper is a generalization of Fomin's
growth diagram ideas to the $K$-theory context.

Let ${\mathbb Y}$ be the 
{\bf Young lattice}, the partial order $\preceq$ on all shapes 
where $\lambda\preceq \mu$ when $\lambda$ is contained inside $\mu$. The 
covering relations on ${\mathbb Y}$ are $\lambda\preceq \mu$ such 
that $\mu/\lambda$ is a single box. 

Each increasing tableau $T$ can be viewed as a {\bf shape sequence} 
of increasing
shapes in ${\mathbb Y}$ where each successive shape is grown from the previous
one by adding some number of boxes, no two in the same row or column. 

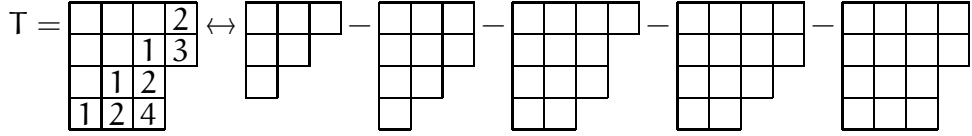
\begin{Example}
\label{exa:shapeseq}
\[T=\tableau{{\ }&{\ }&{\ }&{2}\\{\ }&{\ }&{1}&{3}\\{\ }&{1}&{2}\\{1}&{2}&{4}}
\leftrightarrow
\tableau{{\ }&{\ }&{\ }\\{\ }&{\ }\\{\ }}-
\tableau{{\ }&{\ }&{\ }\\{\ }&{\ }&{\ }\\{\ }&{\ }\\{\ }}-
\tableau{{\ }&{\ }&{\ }&{\ }\\{\ }&{\ }&{\ }\\{\ }&{\ }&{\ }\\{\ }&{\ }}-
\tableau{{\ }&{\ }&{\ }&{\ }\\{\ }&{\ }&{\ }&{\ }\\{\ }&{\ }&{\ }\\{\ }&{\ }}-
\tableau{{\ }&{\ }&{\ }&{\ }\\{\ }&{\ }&{\ }&{\ }\\{\ }&{\ }&{\ }\\{\ }&{\ }&{\ }}\]
\end{Example}

Now, consider the following choice of rectification order:
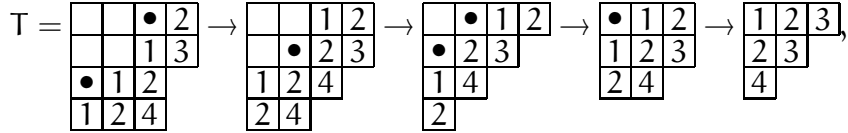
\[
T=\tableau{{\ }&{\ }&{\bullet }&{2}\\{\ }&{\ }&{1}&{3}\\{\bullet }&{1}&{2}\\{1}&{2}&{4}}\to 
\tableau{{\ }&{\ }&{1}&{2}\\{\ }&{\bullet }&{2}&{3}\\{1}&{2}&{4}\\{2}&{4}}\to
\tableau{{\ }&{\bullet }&{1}&{2}\\{\bullet }&{2}&{3}\\{1}&{4}\\{2}}\to
\tableau{{\bullet }&{1}&{2}\\{1}&{2}&{3}\\{2}&{4}}\to
\tableau{{1}&{2}&{3}\\{2}&{3}\\{4}},
\]
where the $\bullet$'s indicate the set of boxes to use in each $K{\tt jdt}$ step.
Each of these increasing tableaux also has a shape sequence, 
which we put one 
atop of another so the shapes increase moving up and to the right.
The result is a {\bf $K$-theory growth diagram}; in our example, we have:

\begin{table}[h]
\begin{center}
\begin{tabular}{|l|l|l|l|l|}
\hline
$(3,2,1)$ & $(3,3,2,1)$ & $(4,3,3,2)$ & $(4,4,3,2)$ & $(4,4,3,3)$\\ \hline
$(2,2)$ & $(3,2,1)$ & $(4,3,2,1)$ & $(4,4,2,1)$ & $(4,4,3,2)$\\ \hline
$(2,1)$ & $(3,1,1)$ & $(4,2,1,1)$ & $(4,3,1,1)$ & $(4,3,2,1)$\\ \hline
$(1)$ & $(2,1)$ & $(3,2,1)$ & $(3,3,1)$ & $(3,3,2)$ \\ \hline
$\emptyset$ & $(1)$ & $(2,1)$ & $(3,2)$ & $(3,2,1)$\\ \hline
\end{tabular}
\end{center}
\caption{A $K$-theory growth diagram}
\end{table}

Consider the following local conditions on any $2\times 2$
subsquare $\tableau{{\alpha}&{\beta}\\{\gamma}&{\delta}}$ 
of such a grid of shapes, where by assumption 
$\gamma\subseteq\alpha\subseteq\beta$ and 
$\gamma\subseteq\delta\subseteq\beta$, as in the example above:
\begin{itemize}
\item[(G1)] $\alpha/\gamma$ is a collection of boxes no two in the same
row or column, and similarly for $\beta/\alpha$, $\beta/\delta$, and
$\delta/\gamma$.  
\item[(G2)] $\delta=K{\tt jdt}_{\alpha/\gamma}(T)$ 
where $T$ is the filling of $\beta/\alpha$ by 1's.
This uniquely determines $\delta$ from $\gamma,\alpha$ and
$\beta$. Similarly, $\alpha$ is uniquely determined by $\gamma, \delta$ and 
$\beta$.
\end{itemize}
In particular, (G1) and (G2) are symmetric in $\alpha$ and 
$\delta$.

\begin{Proposition}\label{prop:gr}
If $\tableau{{\alpha}&{\beta}\\{\gamma}&{\delta}}$ is a $2\!\times\! 2$
square in a $K$-theory growth diagram, then (G1) and (G2) hold.
\end{Proposition}
\begin{proof}
This is a straightforward verification.
\end{proof}

Note therefore that if $\mathcal G$ is a growth diagram, then so is 
$\mathcal G$ reflected about its antidiagonal.

Let $K{\tt GROWTH}(\lambda,\mu;\nu)$ be the number of $K$-theory
growth diagrams such that:
\begin{itemize}
\item the leftmost column
encodes the superstandard tableau of shape $\lambda$;
\item the bottom-most row encodes the superstandard tableau of shape $\mu$;
\item the top right corner is the shape $\nu$.
\end{itemize}

The following fact is immediate from Theorem~\ref{thm:KLR}, giving an alternative formulation of Theorem~\ref{thm:KLR}:

\begin{Corollary}
\label{cor:Kgrowthrule}
(of Theorem~\ref{thm:KLR})
$(-1)^{|\nu|-|\lambda|-|\mu|}C_{\lambda,\mu}^{\nu}
=\#K{\tt GROWTH}(\lambda,\mu;\nu)$.
\end{Corollary}

By the symmetry of growth diagrams, the roles of the $\lambda$ and $\mu$ 
can be interchanged,
resulting in the same growth diagram (up to reflection). Therefore, 
the
rule of Corollary~\ref{cor:Kgrowthrule}
manifests the ${\mathbb Z}_2$ {\bf commutation 
symmetry} $C_{\lambda,\mu}^{\nu}=C_{\mu,\lambda}^{\nu}$
coming from  $[{\mathcal O}_{X_{\lambda}}][{\mathcal O}_{X_{\mu}}]=
[{\mathcal O}_{X_{\mu}}][{\mathcal O}_{X_{\lambda}}]$.

An ${\mathbb Z}_3$-symmetric
rule preserving the {\bf triality symmetry} 
\[C_{\lambda,\mu,\nu^{\vee}}=C_{\mu,\nu^{\vee},\lambda}=C_{\nu^{\vee},\lambda,\mu}\]
where $C_{\lambda,\mu,\nu^{\vee}}:=C_{\lambda,\mu}^{\nu}$ etc., exists in
the form of \emph{puzzles}; see \cite{Vakil}). (Unlike in cohomology, in $K$-theory,
this latter symmetry is not immediate from the geometric definitions; for a proof see
\cite{Buch:KLR, Vakil}. In fact, this symmetry is not expected to hold for 
general $G/P$, although
A.~Knutson has informed us, in private communication, 
that it holds in the minuscule setting.)
One can also hope for a manifestly $S_3$-symmetric rule, 
as is available for the cohomological 
Littlewood-Richardson coefficients via \emph{cartons} 
\cite{Thomas.Yong:IV}. 
However we do not know how to extend cartons to the present context; 
see Section~8 for more notes.

Growth diagrams corresponding to the classical 
rectifications of a standard tableau (using only ${\tt jdt}$ moves) were
first introduced by Fomin, see \cite[Appendix 1]{Stanley} and the
references therein. In that case, the above Proposition simplifies. Specifically,
\begin{itemize}
\item[(F1)] shapes increase by precisely one box in the ``up'' and ``right'' directions.
\item[(F2)] if $\alpha$ is the unique shape containing $\gamma$ and
contained in $\beta$, then $\delta=\alpha$;
otherwise there is a unique such shape different than $\alpha$, 
and this shape is $\delta$. 
\end{itemize}
(Similarly, $\alpha$ is uniquely determined by $\beta, \gamma$ and $\delta$.)

Fomin's growth diagrams provide further useful combinatorial ideas 
that we extend below to the $K$-theory setting. These diagrams 
also arise (along with other classical tableaux algorithms we generalize) 
in an elegant geometric context, 
due to work of van Leeuwen \cite{Leeuwen}; 
there are reasons to hope that one can 
extend his work to the setting of this paper.

\section{The infusion involution}

Given $T\in {\tt INC}(\lambda/\alpha)$ and $U\in {\tt INC}(\nu/\lambda)$
define 
\[K{\tt infusion}(T,U)=(K{\tt infusion}_{1}(T,U), K{\tt infusion}_2(T,U))
\in {\tt INC}(\gamma/\alpha)\times {\tt INC}(\nu/\gamma)\]
(for some straight shape $\gamma$) 
as follows: consider the largest label ``$m$'' that appears in $T$,
appearing at $x_1,\ldots,x_k$.  Apply the slide $K{\tt jdt}_{\{x_i\}}(U)$,
leaving some ``holes'' at the other side of $\nu/\lambda$.
Place ``$m$'' in these holes and repeat, moving the labels
originally from $U$ until all labels of $T$
are exhausted. The resulting
tableau of shape $\gamma/\alpha$ and skew
tableau of shape $\nu/\gamma$ are the outputted tableaux.
To define 
\[K{\tt revinfusion}(T,U)\!=\!(K{\tt revinfusion}_1(T,U), \!K{\tt revinfusion}_{2}(T,U))
\!\in\!\! {\tt INC}(\gamma/\alpha)\!\times\! {\tt INC}(\nu/\gamma)\]
we apply $K{\tt revjdt}$ moves to $T$, moving
into boxes of $U$. We begin by removing the labels ``$1$'' appearing in $U$ at
boxes $\{x_i\}\in \nu/\lambda$, apply ${\tt revjdt}_{\{x_i\}}(T)$, 
and place the ``$1$'' in the vacated holes of $\lambda$ and 
continuing with higher labels of $U$.

It is easy to show
$K${\tt infusion} and $K${\tt revinfusion} are inverses of one
another, by inductively applying the observation that if $\{y_i\}$ are the
boxes vacated by $K{\tt jdt}_{\{x_i\}}(T)$ then
\[K{\tt revjdt}_{\{y_i\}}(K{\tt jdt}_{\{x_i\}}(T))=T.\]
We will need the following fact (the ``infusion involution'', cf.~\cite{Haiman:DE, Sottile:tabswitch}) : 
\begin{Theorem}
\label{thm:infusionandrev}
For any increasing tableaux $T$ and $U$ such that ${\tt shape}(U)$ extends
${\tt shape}(T)$ then $K{\tt infusion}(T,U)=K{\tt revinfusion}(T,U)$. That is,
$K{\tt infusion}(K{\tt infusion}(T,U))=(T,U)$.
\end{Theorem}

\begin{Example}
If $T=\tableau{{\underline 1}&{\underline 2}&{\underline 3}\\{\underline 2}&
{\underline 3}\\{\underline 4}}$ and
$U=\tableau{&&&{2}\\&&{1}&{3}\\&{1}&{3}\\{2}&{3}&{4}}$ then we compute
$K${\tt infusion} as follows:

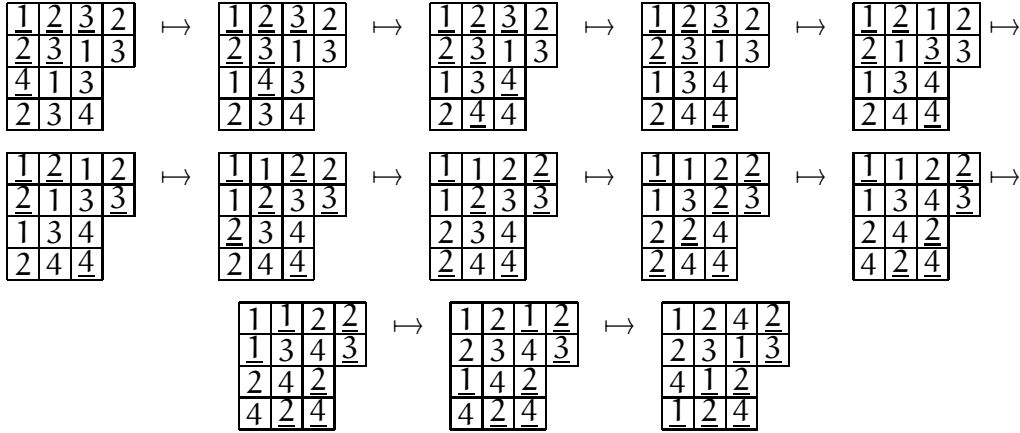
\[\begin{array}{ccccccccc}\nonumber
\tableau{{\underline 1}&{\underline 2}&{\underline 3}&{2}\\{\underline 2}&{\underline 3}&{1}&{3}\\
{\underline 4}&{1}&{3}\\{2}&{3}&{4}}
& \mapsto & 
\tableau{{\underline 1}&{\underline 2}&{\underline 3}&{2}\\{\underline 2}&{\underline 3}&{1}&{3}\\{1}&{\underline 4}&{3}\\{2}&{3}&{4}}
& \mapsto & 
\tableau{{\underline 1}&{\underline 2}&{\underline 3}&{2}\\{\underline 2}&{\underline 3}&{1}&{3}\\{1}&{3}&{\underline 4}\\{2}&{\underline 4}&{4}}
& \mapsto & 
\tableau{{\underline 1}&{\underline 2}&{\underline 3}&{2}\\{\underline 2}&{\underline 3}&{1}&{3}\\{1}&{3}&{4}\\{2}&{4}&{\underline 4}}
& \mapsto & 
\tableau{{\underline 1}&{\underline 2}&{1}&{2}\\{\underline 2}&{1}&{\underline 3}&{3}\\{1}&{3}&{4}\\{2}&{4}&{\underline 4}}\mapsto\\
\end{array}
\]
\[
\begin{array}{ccccccccc}\nonumber
\tableau{{\underline 1}&{\underline 2}&{1}&{2}\\{\underline 2}&{1}&{3}&{\underline 3}\\
{1}&{3}&{4}\\{2}&{4}&{\underline 4}}
& \mapsto & 
\tableau{{\underline 1}&{1}&{\underline 2}&{2}\\{1}&{\underline 2}&{3}&{\underline 3}\\{\underline 2}&{3}&{4}\\{2}&{4}&{\underline 4}}
& \mapsto & 
\tableau{{\underline 1}&{1}&{2}&{\underline 2}\\{1}&{\underline 2}&{3}&{\underline 3}\\{2}&{3}&{4}\\{\underline 2}&{4}&{\underline 4}}
& \mapsto & 
\tableau{{\underline 1}&{1}&{2}&{\underline 2}\\{1}&{3}&{\underline 2}&{\underline 3}\\{2}&{\underline 2}&{4}\\{\underline 2}&{4}&{\underline 4}}
& \mapsto & 
\tableau{{\underline 1}&{1}&{2}&{\underline 2}\\{1}&{3}&{4}&{\underline 3}\\{2}&{4}&{\underline 2}\\{4}&{\underline 2}&{\underline 4}}\mapsto\\
\end{array}\] 
\[
\begin{array}{ccccccccc}\nonumber
\tableau{{1}&{\underline 1}&{2}&{\underline 2}\\{\underline 1}&{3}&{4}&{\underline 3}\\
{2}&{4}&{\underline 2}\\{4}&{\underline 2}&{\underline 4}}
& \mapsto & 
\tableau{{1}&{2}&{\underline 1}&{\underline 2}\\{2}&{3}&{4}&{\underline 3}\\
{\underline 1}&{4}&{\underline 2}\\{4}&{\underline 2}&{\underline 4}}
& \mapsto & 
\tableau{{1}&{2}&{4}&{\underline 2}\\{2}&{3}&{\underline 1}&{\underline 3}\\
{4}&{\underline 1}&{\underline 2}\\{\underline 1}&{\underline 2}&{\underline 4}}\\
\end{array}\] 
Hence 
\[K{\tt infusion}(T,U)=\left(\tableau{{1}&{2}&{4}\\{2}&{3}\\{4}},
\tableau{&&&{\underline 2}\\&&{\underline 1}&{\underline 3}\\&{\underline 1}&{\underline 2}\\{\underline 1}&{\underline 2}&{\underline 4}}\right).\]
The reader can check that applying $K{\tt infusion}$ to this pair returns
$(T,U)$, in agreement with the Theorem.
\end{Example}

\begin{proof}
Construct the growth diagram for ${K}{\tt rect}(U)$ using the
slides suggested by the entries of $T$. Notice the bottom row represents
$K{\tt infusion}_{1}(T,U)$ and the right column
represents $K{\tt infusion}_{2}(T,U)$. However, by the antidiagonal
symmetry of growth diagrams, the growth diagram computing 
$K{\tt infusion}$ applied to 
$K{\tt infusion}(T,U)$ is simply the one for $K{\tt infusion}(T,U)$ 
reflected about the antidiagonal.
\end{proof}

Finally, the growth diagram formalism makes it straightforward to observe
facts such as the following, which we will need in Section~6:

\begin{Lemma} 
\label{lemma:split}
Let $T\in {\tt INC}(\nu/\lambda)$, $R\in{\tt INC}(\lambda)$ and  
fix $a\in {\mathbb N}$. If $A$ be the increasing tableau 
consisting of entries
from $1$ to $a$ of $T$, and $B=T\setminus A$ is the remaining tableau,  
then $K{\tt
infusion}_1(R,T)=K{\tt infusion}_1(R,A)\cup
K{\tt infusion}_1(K{\tt infusion}_2(R,A),B)$.
\end{Lemma}
\begin{proof}
Draw the growth diagram for
$K{\tt infusion}(R,T)$, encoding $R$ on the left and $T$ on the top.
The shape ${\tt shape}(R)\cup{\tt shape}(A)$ appears on the top row.
Draw a vertical line through the growth diagram at that point.  The
diagram to the left of this line encodes the rectification of $A$ by
$R$.  The diagram to the right of the line encodes the rectification
of $B=T\setminus A$ by the tableau encoded along the dividing line, which
is $K{\tt infusion}_2(R,A)$.
\end{proof}

\section{A generalization of Sch\"{u}tzenberger's evacuation involution}

While on the topic of growth diagrams, we 
take this opportunity to introduce a generalization of another
classical result from tableau theory. This section will not be needed
in the remainder of the paper.

For $T\in {\tt INC}(\lambda)$, let 
$^{\circ}T$ be obtained by erasing the (unique) entry $1$ 
in the northwest corner $c$ of $T$ and subtracting $1$ from the 
remaining entries. Let
\[\Delta(T)=K{\tt jdt}_{\{c\}}({^{\circ}T}).\]
The $K$-{\bf evacuation} $K{\tt evac}(T)\in {\tt INC}(\lambda)$ 
is defined by the shape sequence
\[\emptyset={\tt shape}(\Delta^{\max T}(T))-
{\tt shape}(\Delta^{\max T-1}(T))-\ldots - {\tt shape}(\Delta^{1}(T)) - T.\]

The following result extends Sch\"{u}tzenberger's classical theorem
for $T\in {\tt SYT}(\lambda)$.

\begin{Theorem}
\label{thm:evac}
$K{\tt evac}:{\tt INC}(\lambda)\to {\tt INC}(\lambda)$
is an involution, i.e., $K{\tt evac}(K{\tt evac}(T))=T$.
\end{Theorem}

\begin{Example}
\label{exa:schutgrowth}
Let $T=\tableau{{1}&{2}&{3}&{5}\\{2}&{3}&{4}\\{4}&{5}}\in 
{\tt INC}((4,3,2))$. Then the $K$-evacuation is computed by 
\[
\Delta^{1}(T)=\tableau{{1}&{2}&{3}&{4}\\{2}&{3}\\{3}&{4}}\mapsto
\Delta^{2}(T)=
\tableau{{1}&{2}&{3}\\{2}&{3}\\{3}}\mapsto
\Delta^{3}(T)=\tableau{{1}&{2}\\{2}}\mapsto
\Delta^{4}(T)=\tableau{{1}}\mapsto \Delta^{5}(T)\mapsto \emptyset.
\]
Thus $K{\tt evac}(T)=\tableau{{1}&{2}&{3}&{4}\\{2}&{3}&{5}\\{3}&{4}}$.
One checks that applying $K{\tt evac}$ to this tableau returns $T$.
\end{Example}

\noindent
\emph{Proof of Theorem~\ref{thm:evac}:}
Express each of the increasing tableaux
\[T,\Delta^{1}(T),\ldots, \Delta^{\max T-1}(T),\Delta^{\max T}(T)=\emptyset\]
as a shape sequence and place them right justified in a triangular growth
diagram. In the example above, we have Table~\ref{table:triangle}.
Noting that each ``minor'' of the table whose southwest corner contains a
``$\emptyset$'' is in fact a growth diagram, it follows that the triangular
growth diagram can be reconstructed using (G1) and (G2), by
Proposition~\ref{prop:gr}. 
Observe that the right column encodes
$K{\tt evac}(T)$. By the symmetry of growth diagrams,
it follows that applying the above procedure to $K{\tt evac}(T)$ 
would give the same triangular growth diagram, after a reflection
across the antidiagonal. Thus the result follows.\qed

\begin{table}[h]
\begin{center}
\begin{tabular}{|l|l|l|l|l|l|}
\hline
$\emptyset$ & $(1)$ & $(2,1)$ & $(3,2)$ & $(3,3,1)$ & $(4,3,2)$\\ \hline
& $\emptyset$ & $(1)$ & $(2,1)$ & $(3,2,1)$ & $(4,2,2)$\\ \hline
& & $\emptyset$ & $(1)$ & $(2,1)$ & $(3,2,1)$\\ \hline
& & & $\emptyset$ & $(1)$ & $(2,1)$ \\ \hline
& & & & $\emptyset$ & $(1)$\\ \hline
& & & & & $\emptyset$ \\ \hline
\end{tabular}
\end{center}
\caption{\label{table:triangle} A triangular growth diagram for Example~\ref{exa:schutgrowth}.}
\end{table}

\section{Proof of the $K{\tt jdt}$ rule}

The strategy of our proof is based on the following fact.
In the cohomological context, this approach was utilized 
in~\cite{KTW:EJC, BKT}.

\begin{Lemma}
Let $\{d_{\lambda,\mu}^{\nu}\}$ be integers indexed by shapes
$\lambda,\mu,\nu\subseteq \Lambda$ that: 
\begin{itemize} 
\item[(A)] define a commutative and associative
ring $(R,\circ)$, with ${\mathbb Z}$-basis $\{a_{\lambda}\}$ 
indexed by shapes $\lambda\subseteq \Lambda$, by: 
\[a_{\lambda}\circ a_{\mu}=\sum_{\nu\subseteq \Lambda} d_{\lambda,\mu}^{\nu}a_{\nu}, \mbox{ and }\] 
\item[(B)] $d_{\lambda,\rho}^{\nu}=c_{\lambda,\rho}^{\nu}$ whenever $\rho=(t)$ for $0\leq t\leq n-k$
\end{itemize}
then $d_{\lambda,\mu}^{\nu}=c_{\lambda,\mu}^{\nu}$. 
\end{Lemma}
\begin{proof}
The class $[\mathcal O_{X_\lambda}]$ can be expressed as a polynomial
in $[\mathcal O_{X_{(1)}}],\dots,[\mathcal O_{X_{(n-k)}}]$. 
This follows
by an easy downward induction on $|\lambda|$ using the fact that  
such an expression
exists in cohomology for $[X_\lambda]\in H^{\star}(X,{\mathbb Q})$ as 
a polynomial in the classes $[X_{(t)}]$ (the \emph{Jacobi-Trudi identity}) and
the lowest order term in 
$K$-theory agrees with cohomology under the Chern isomorphism. 
Let this
polynomial be $P_{\lambda}(X_1,\ldots,X_{n-k})$ (where above 
$X_t=[{\mathcal O}_{X_{(t)}}]$).
Now (A) and (B) imply $a_{\lambda}=P_{\lambda}(a_{(1)},\ldots, a_{(t)})$.
Using (B) again, we see that 
the map from $(R,\circ)$ to $K(X)$ sending $a_{\lambda}\mapsto
[{\mathcal O}_{X_{\lambda}}]$ is a ring isomorphism, so the desired 
conclusion follows.
\end{proof}

To apply the lemma, let $d_{\lambda,\mu}^{\nu}$ be the integers computed
by the rule given in the statement of the theorem. It remains to check
associativity and agreement with Pieri's rule, which we do below.
Note,  in our proof of associativity we \emph{assume}
Theorem~\ref{thm:welldefined} is true -- this latter result 
is actually proved in the following section, using some of 
the elements introduced in the proof of agreement of Pieri's rule, which
of course, does not use this assumption.

\noindent
\emph{Associativity:} 
Let $\alpha$, $\beta$, $\gamma$, $\nu$ be straight shapes and
fix superstandard tableaux $S_{\alpha}, S_{\beta},
S_{\gamma}$ and $S_\nu$.
 
Associativity is the assertion that
\begin{equation}
\label{eqn:assoc}
\sum_\sigma d_{\alpha,\beta}^\sigma d_{\sigma,\gamma}^\nu=
\sum_\tau d_{\alpha,\tau}^\nu d_{\beta,\gamma}^\tau.
\end{equation}

The lefthand side of (\ref{eqn:assoc}) 
counts pairs of tableaux $(B,C)$ where $B$ is of shape $\sigma/\alpha$
such that $K{\tt rect}(B)=S_\beta$, 
and $C$ is of shape $\nu/\sigma$ such that
$K{\tt rect}(C)=S_\gamma$.

Let $K{\tt infusion}(S_{\alpha},B)=(S_{\beta},A)$ where $A$ is of
shape $\sigma/\beta$, and $K{\tt rect}(A)=S_\alpha$. Next compute
$K{\tt infusion}(A,C)=(D,E)$. We have that
$K{\tt rect}(E)=S_\alpha$ (since this was the case with $A$) and that
${\tt shape}(E)=\nu/\tau$ for some $\tau$, and similarly 
$K{\tt rect}(D)=S_\gamma$ (since this was the case for $C$) and
${\tt shape}(D)=\tau/\beta$.

By Theorem~\ref{thm:infusionandrev} it follows that
the above process establishes a bijection 
\[(B,C)\mapsto (E,D)\]
into the set of pairs of tableaux counted by the righthand side of
(\ref{eqn:assoc}).  (More precisely, for pairs counted by
\[\sum_\tau d_{\tau,\alpha}^\nu d_{\beta,\gamma}^\tau=
\sum_\tau d_{\alpha,\tau}^\nu d_{\beta,\gamma}^\tau\]
where the equality $d_{\tau,\alpha}^\nu=d_{\alpha,\tau}^\nu$
follows from an easy argument using $K{\tt infusion}$.)
Associativity follows.  

\noindent
\emph{Agreement with Pieri's rule:} We prove our rule agrees with the 
following formula, due to Lenart~\cite{Lenart}:

\begin{Theorem}
\label{thm:lenart}
Let $r(\nu/\lambda)$ be the number of rows of $\nu/\lambda$. Then
\[[{\mathcal O}_{X_{\lambda}}][{\mathcal O}_{X_{(t)}}]=
\sum_{\nu}(-1)^{|\nu|-|\mu|-t}{r(\nu/\lambda)-1 \choose |\nu/\lambda|-t}[{\mathcal O}_{X_{\nu}}],\]
where the sum ranges over all $\nu\subseteq \Lambda$ obtained by
adding a horizontal strip (no two added boxes are in the same column)
to $\lambda$ of size at least $t$.
\end{Theorem}

Our task is to show that $d_{\lambda,(t)}^{\nu}={r(\nu/\lambda)-1 \choose |\nu/\lambda|-t}$ when $\nu$ is of the form in the statement of 
Theorem~\ref{thm:lenart} and is zero otherwise.

First assume $\nu$ is of the desired form. Note that if 
$|\nu/\lambda|-t>r(\nu/\lambda)-1$ then no increasing filling by
$\{1,\dots,t\}$ is possible,
as desired. So assuming otherwise, we proceed to construct the required
number of increasing tableaux on $\nu/\lambda$, as follows. Select 
$|\nu/\lambda|-t$ of the non-bottom-most
$r(\nu/\lambda)-1$ rows of $\nu/\lambda$. Now fill the bottom row with
consecutive entries $1,2,\ldots,k$ where $k$
is the number of boxes in that bottom row of $\nu/\lambda$. Proceeding to fill the 
remaining boxes of $\nu/\lambda$ from southwest to northeast. If the current
row to be filled was one of the $|\nu/\lambda|-t$ selected rows then
begin with the last entry $e$ used in the previously filled row. Otherwise
use $e+1$.

Call thse fillings $t$-{\bf Pieri fillings}.

\begin{Example}
Suppose $\lambda=(5,3,2)$, $\nu=(6,5,2,2)$ and $t=4$. Then $r(\nu/\lambda)=3$
and $|\nu/\lambda|-t=1$. Hence the two $4$-Pieri fillings we construct are
\[\tableau{{ \ }&{\ }&{\ }&{\ }&{\ }&{4}\\{\ }&{\ }&{\ }&{2}&{3}\\
{\ }&{\ }\\{1}&{2}} \mbox{ \ \ \ and\ \ \  } 
\tableau{{ \ }&{\ }&{\ }&{\ }&{\ }&{4}\\{\ }&{\ }&{\ }&{3}&{4}\\
{\ }&{\ }\\{1}&{2}},
\]
which both rectify to $\tableau{{1}&{2}&{3}&{4}}$. (In the first
tableau we selected the second row and in the second we selected the
top row.) 
\end{Example}

\begin{Lemma}
\label{lemma:tPieri}
For any rectification order, 
a $t$-Pieri filling $K$-rectifies to $S_{(t)}$.
No other increasing tableau $K$-rectifies to $S_{(t)}$, for any choice of
rectification order.
\end{Lemma}
\begin{proof}
That the $t$-Pieri fillings all $K$-rectify (under any rectification
order) to $S_{(t)}$ follows from a 
straightforward induction on $|\lambda|\geq 0$ where we show in fact
that any $K{\tt jdt}$ slide applied to a $t$-Pieri filling results in
a $t$-Pieri filling.

A similar induction
shows that no other increasing tableau from ${\tt INC}(\nu/\lambda)$ 
$K$-rectifies to $S_{(t)}$ (noting that any such tableau with entries in
$\{1,\dots,t\}$ has a pair of entries
$i<j$ where $j$ is southwest of $i$). Separately, but for similar reasons,
when $\nu/\lambda$ is not
a horizontal strip, one more induction on $|\lambda|$ proves
no increasing tableau can $K$-rectify to $S_{(t)}$.
\end{proof}

This completes the proof of Theorem~\ref{thm:KLR}, assuming Theorem~\ref{thm:welldefined}.\qed

\section{Proof of the $K{\tt rect}$ theorem}

We now prove Theorem~\ref{thm:welldefined}. First define the 
{\bf reading word} of a tableau $T$ to be the word obtained by reading the
rows of $T$ from left to right, starting from the bottom and moving up.
Let ${\tt LIS}(T)$ be the length of the longest strictly increasing 
subsequence of the reading word of $T$.  

The following result is crucial to our proof of Theorem~\ref{thm:welldefined}.

\begin{Theorem}
\label{thm:LIS}
${\tt LIS}({K}{\tt jdt}_{\{x_i\}}(T))={\tt LIS}(T)$. In particular,
any rectification order applied to $T$ results in a straight shape whose
first row has length equal to ${\tt LIS}(T)$.
\end{Theorem}

\begin{Example}
Consider the two (different) rectifications of the same tableau $T$
performed in Example~\ref{exa:twodiff}. The reading word of $T$ is
${\underline 1} \ {\underline 3} \ {\underline 4} \ 2 \ 2$ (where the
unique longest strictly increasins subsequence has been underlined) so 
${\tt LIS}(T)=3$. Note that also
${\tt LIS}(T_1)={\tt LIS}(T_2)=3$, that is the lengths of the first rows
of $T_1$ and $T_2$ agree, although $T_1\neq T_2$.
\end{Example}

\noindent
\emph{Proof of Theorem~\ref{thm:LIS}:}
We will show that if $I$ is a set of boxes of $T$ which forms a strictly
increasing subsequence of the reading word of $T$, then there is a string of
boxes of equal length in $K{\tt jdt}_{\{x_i\}}(T)$ which also forms
a strictly increasing subsequence of the reading word.  
A symmetric argument using reverse slides gives the other desired inequality,
thereby implying the theorem.

Fix $I$ as above.
We will analyze the slide ${K}{\tt jdt}_{\{x_i\}}(T)$,
{\tt switch} by {\tt switch}.
Set $T_0:=T$, 
and let
$T_i$ be the result of switching the $\bullet$'s and the $i$'s of $T_{i-1}$.
Initially set $I_0:=I$. 
In a moment, we will describe $I_i$ as 
a collection of some of the boxes of $T_i$.  

We will show that, at each step, $I_i$ has the following properties:
\begin{itemize}
\item[(P1)]
The labels of $I_i$ are strictly 
increasing in the reading word order, except for
perhaps one $\bullet$ box.  
\item[(P2)] If $I_i$ contains a $\bullet$ box, then 
the labels in $I_i$ preceding the $\bullet$ box in the reading word order
are weakly less than $i$, while the labels of boxes following the $\bullet$ box 
are strictly greater than $i$.  
\item[(P3)] If there is a $\bullet$ box $y_i$ in $I_i$, then 
there must be some box $z_i$ in $I_i$, in the same row as $y_i$ and weakly
to its right, such that
the entry in the box  $a_i$ immediately below $z_i$ is strictly
smaller than the entry in the 
next box $b_i$ of $I_i$ after $z_i$, in the reading word order.  
\end{itemize}

\begin{Example}
(P1) and (P2) are self explanatory. For (P3), a possible
configuration that can arise in our discussion below is 
\[\tableau{{\underline 1}&{\underline\bullet}&{\underline 2}&{4}&{\underline 5}&{7}&{\underline 9}\\{\bullet}&{2}&{3}&{6}&{8}&{9}},\] 
where the underlined labels
indicate members of $I_1$. Here the role of $z_1$ is played by the 
$\underline 5$, so $a_1$ is the $8$ and $b_i$ is the ${\underline 9}$. Note
that $b_i$ need not be immediately to the right of the $z_i$.
\end{Example}

We need the following:

\begin{Lemma}\label{lemma:adj} If $I_{i-1}$ satisfies (P1)--(P3) and contains 
an $\bullet$ box and
a box labelled $i$ 
then the $i$ is immediately to the right of the $\bullet$ box.  
\end{Lemma}

\begin{proof} First we use (P3): since the box $a_{i-1}$
below $z_{i-1}$ is southeast of the $\bullet$ box $y_{i-1}$, 
its label cannot be strictly less than $i$, because we are in $T_{i-1}$.
Thus the label of $b_{i-1}$ is greater than $i$.  
By (P1) and (P2) combined, this
implies that the $i$ in $I_{i-1}$ 
must be in the same row as the $\bullet$ box, and
therefore must be immediately to its right (again, for the reason that 
we are in $T_{i-1}$).  
\end{proof}

We now proceed to define $I_i$ inductively for $i\geq 1$.  
Assume that $I_{i-1}$ 
satisfies (P1)--(P3).
After performing the slide interchanging $\bullet$ boxes with $i$'s we define
$I_i$ as follows:

\begin{itemize}
\item[(i)] If $I_{i-1}$ has no box containing $i$, 
then $I_i:=I_{i-1}$.

\item [(ii)] If $I_{i-1}$ has a box containing $i$ and a $\bullet$ box
(i.e., we are in the position of Lemma~\ref{lemma:adj}), then 
$I_i:=I_{i-1}$.  

\item[(iii)] If $I_{i-1}$ has a box containing $i$, but does not
have a $\bullet$ box,
and the $i$ in $I_{i-1}$ does not move, then $I_i:=I_{i-1}$.

\item[(iv)] If $I_{i-1}$ has a box containing $i$, but does not have a
$\bullet$ box,
and there is a $\bullet$ box (not in $I_{i-1}$) 
immediately to the left of the $i$ in 
$I_{i-1}$, then
let $I_i$ be $I_{i-1}$ with the box containing $i$ in $I_{i-1}$ 
replaced by the box to its left (into which $i$ has moved).  

\item[(v)] If $I_{i-1}$ has a box containing $i$, but does not have a
$\bullet$ box,
there is a $\bullet$ box (not in $I_{i-1}$) 
immediately above the $i$, and we are not in case (iv), then let
$I_i$ be $I_{i-1}$ with the box containing $i$ in $I_{i-1}$ and
all the other boxes in $I_{i-1}$ to the right of it in the same row, 
replaced by
the boxes immediately above them.  
\end{itemize}

Clearly (i)-(v) indeed enumerates all of the 
intermediate possibilities during a $K{\tt jdt}$ slide.

We now prove that $I_i$ satisfies (P1)--(P3).

Case (i):
We split this case up into three subcases.  First, we consider the
case that $I_{i-1}$ has no $\bullet$ box.  In this case, (P1) is 
trivially satisfied (since it held for $I_{i-1}$), 
and (P2) and (P3) are vacuously true.  

Next, we consider the subcase that $I_{i-1}$ has
 a $\bullet$ box into which an $i$ (not in $I_{i-1}$) moves in.  
Since (P1) and (P2) are satisfied for $I_{i-1}$, (P1) will be 
satisfied after this, and (P2) and (P3) are vacuous since 
$I_i$ has no $\bullet$ box.  

Finally, we consider the subcase where $I_{i-1}$ has a $\bullet$ box which
stays as such in $I_i$.  Since the contents of $I_{i-1}$ and $I_i$ are
the same,  (P1) and (P2) are satisfied.  To show (P3) is satisfied, 
observe that the label in the box below $z_{i-1}$ is strictly greater 
than $i$ (otherwise $z_{i-1}$ has a label weakly smaller than $i-1$ and is
southeast of a $\bullet$, a contradiction),
so it does not move, and thus we can take $z_i:=z_{i-1}$.  

Using Lemma~\ref{lemma:adj}, it is clear that case (ii) preserves (P1) and (P2). 
To check (P3), as in the previous case, we can take $z_i:=z_{i-1}$.
This would not work if $z_{i-1}=y_{i-1}$, but this is impossible, because
the entry in the box below $z_{i-1}$ should be less than the next entry in
$I_{i-1}$ after $z_{i-1}$, which is $i$.  So the $\bullet$ box is 
immediately above a box which is at most $i-1$, and this can't happen
in $T_{i-1}$.  

Cases (iii) and (iv) are trivial: (P1) holds since the contents of  
$I_{i-1}$ and $I_i$ are the same, and (P2) and (P3) are vacuously true
since $I_i$ contains no $\bullet$ box.

Now we consider case (v).  (P1) is trivial, so if $I_i$ has no
$\bullet$ box, then we are done.  So assume it does.  The only way a
$\bullet$ box could appear in $I_i$ is in the following situation:

\[\tableau{{\bullet }&{i}\\{i}&{j}}\mapsto 
\tableau{{i }&{\bullet}\\{\bullet}&{j}},\]
where the box containing $j$ is also in $I_{i-1}$.

In this situation the the top two boxes will be in $I_i$, and so we will
have introduced a $\bullet$ box into $I_i$. 
(P2) is clearly satisfied. 
Set $z_i$ to be the rightmost of the
boxes that are in $I_i$ but not in $I_{i-1}$.  Now (P3) is satisfied 
because (P1) was satisfied for $I_{i-1}$.  

This completes the proof that $I_i$ satisfies (P1)--(P3).  Thus
after iteration, we eventually terminate
with a set of boxes $I_m$ in $T_m:=K{\tt jdt}_{\{x_i\}}(T)$ which
satisfies (P1)--(P3).  We wish to show that $I_m$ contains no
$\bullet$ box.  Suppose that it did.  This $\bullet$ box of $I_m$
must be an outer corner of $T$ (by the way $K{\tt jdt}$ is defined).  
This contradicts (P3), since the square
below $z_i$ is southeast of the $\bullet$ box, and thus contains no
label.
Thus $I_m$ contains no $\bullet$ box, so (P1) implies that 
there is a strictly increasing subsequence of the reading
word of $K{\tt jdt}_{\{x_i\}}(T)$ whose length equals 
the length of $I$, as desired.  
\qed

\begin{Remark}
Theorem~\ref{thm:LIS} may be regarded as a generalization of the classical
result of Schensted which asserts that the longest
increasing subsequence of a permutation $w=w_1 w_2\ldots w_n$ in the
symmetric group $S_n$ (written in one-line notation) is equal to the first row
of the common shape of the corresponding insertion and recording tableaux 
under the Robinson-Schensted algorithm; see, e.g.,\cite{Stanley}. To see this,
one needs to use the well-known fact that the insertion tableau of $w$
is equal to the (classical) rectification of the ``permutation tableau''
$T_{w}$ of skew shape 
\[(n,n-1,n-2,\ldots,3,2,1)/(n-1,n-2,\ldots,3,2,1),\]
where $w_1$ occupies the southwest-most box, followed by $w_2$ in the box
to its immediate northeast, etc.
In \cite{Thomas.Yong:VI} we further explore this observation, and connect
${K}{\tt rect}$ to the Hecke algorithm of \cite{BKSTY}.
\end{Remark}

Recall the definition of $t$-Pieri filling given in the previous section.

\begin{Lemma}
\label{lemma:verysplit}
If an increasing tableau $T$ rectifies (with respect to any
rectification order) to a tableau $V$ which has
precisely $1,2,\ldots, t$ in the first row and no labels weakly
smaller than $t$ elsewhere, then:
\begin{enumerate}
\item the labels $1,2,\ldots,t$ form a subtableau of $T$ that is a 
$t$-Pieri filling, and 
\item ${\tt LIS}(T)=t$.  
\end{enumerate}
Conversely, if $T$ satisfies (1) and (2), then its $K$-rectification with
respect to any rectification order must have
$1,2,\ldots,t$ as its first row.
\end{Lemma}
\begin{proof}
By Lemma~\ref{lemma:split}, 
$V$ contains the rectification of the subtableau
of $T$ consisting of the entries between 1 and $t$; by results of the
previous section, it follows that these entries must form a $t$-Pieri
filling; this proves that (1) holds.
By Theorem~\ref{thm:LIS}, 
${\tt LIS}(T)={\tt LIS}(V)=t$, proving (2).  

Now suppose $T$ satisfies (1) and (2). By the discussion of $t$-Pieri
fillings of Section~5, and by Theorem~\ref{thm:LIS},
these properties are
preserved under $K{\tt jdt}$ slides.  Thus, any rectification of
$T$ must satisfy these properties.  Any filling of a straight shape 
satisfying these properties must have first row $1,2,\ldots,t$ and no
entries weakly smaller than $t$ elsewhere. 
\end{proof}

\noindent
\emph{Proof of Theorem~\ref{thm:welldefined}:}
Let $R\in {\tt INC}(\lambda)$ encode a rectification where
$K{\tt infusion}_1(R,T)=S_{\mu}$. Let us suppose that the first
row of $S_{\mu}$ is $S_{(t)}$.
By Theorem~\ref{thm:LIS}, ${\tt LIS}(T)=t$.  By 
Lemma~\ref{lemma:verysplit}, the subtableau $P$ of $T$, consisting of
the boxes containing one of the labels $1,2,\ldots,t$, is a $t$-Pieri
filling.  

Suppose $Q\in {\tt INC}(\lambda)$ is
another rectification order. Since the labels of $P$ and $T\setminus P$ 
contain labels weakly smaller than $t$, a strictly larger than $t$, 
respectively, by Lemma~\ref{lemma:split}, 
we can
compute $V:=K{\tt infusion}_1(Q,T)$ in two stages. 
First, by Lemma~\ref{lemma:tPieri}, $K{\tt infusion}_1(Q,P)$ is simply 
$S_{(t)}$, because $P$ is a 
$t$-Pieri filling. Secondly, 
we use $K{\tt infusion}_2(Q,P)$ to (partially) rectify $T\setminus P$.
{\it A priori}, this could contribute extra boxes to first row of $V$ but
since, by Theorem~\ref{thm:LIS}, ${\tt LIS}(V)={\tt LIS}(T)=t$, it does not.  
Thus the rectification of $T$ by $Q$ consists of the row 
$S_{(t)}$ with a rectification of $T\setminus P$ to
a straight shape underneath it.  

Now, by assumption $T\setminus P$ has a (partial) 
rectification to a superstandard
tableaux (using labels starting from $t+1$), namely $S\setminus S_{(t)}$. So
by induction on the number of boxes of the starting shape, we can conclude
that $T\setminus P$ will (partially) rectify to $S\setminus S_{(t)}$ under any rectification
order. Therefore $V=S$, as desired.
\qed

\section{Minuscule Schubert calculus conjectures: example and discussion}

As stated in the introduction, the minuscule $G/P$ are classified into five infinite families of
spaces (Grassmannians, odd/even orthogonal Grassmannians, odd projective
space, even dimensional quadrics) associated to the classical Lie groups,
and two exceptional cases (the Cayley plane and 
$G_{\omega}({\mathbb O}^3,{\mathbb O}^6)$) associated to the Lie types $E_6$ and $E_7$. 
The interested reader may find details compatible
with the notation used here in~\cite{Thomas.Yong}; in particular, there we concretely
describe $\Lambda_{G/P}$ in each of these cases. Thus, for brevity, we 
content ourselves with an example to illustrate our conjecture.

\begin{Example}
\label{exa:E6}
Let $G/P={\mathbb O}{\mathbb P}^2$ be the Cayley plane. Here we have
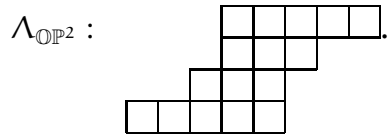
\[\Lambda_{{\mathbb O}{\mathbb P}^2}: \ \ \ \tableau{&&&{\ }&{\ }&{\ }&{\ }&{\ }\\&&&{\ }&{\ }&{\ }\\&&{\ }&{\ }&{\ }\\
{\ }&{\ }&{\ }&{\ }&{\ }}.\]
We conjecturally compute $C_{\lambda,\mu}^{\nu}({\mathbb O}{\mathbb P}^2)$ where
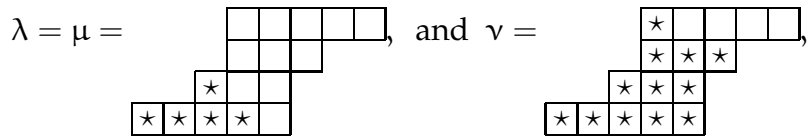
\[\lambda=\mu=\tableau{&&&{\ }&{\ }&{\ }&{\  }&{\ }\\&&&{\ }&{\ }&{\ }\\&&{\star }&{\ }&{\ }\\
{\star }&{\star }&{\star }&{\star }&{\ }}, \mbox{ \ and \ } 
\nu=\tableau{&&&{\star }&{\ }&{\ }&{\  }&{\ }\\&&&{\star }&{\star }&{\star }\\&&{\star }&{\star }&{\star }\\
{\star }&{\star }&{\star }&{\star }&{\star }},
\]
where the southwest-most box is the unique minimum of $\Lambda_{{\mathbb O}{\mathbb P}^2}$ and the poset increases as one moves ``right'' or ``up''.

The relevant shapes/lower order ideals of $\Lambda_{{\mathbb O}{\mathbb P}^2}$ 
are indicated by the boxes filled with
$\star$'s. We can encode the shapes by the size of columns as read from
left to right, so $\lambda=\mu=(1,1,2,1)$ and $\nu=(1,1,2,4,3,1)$. Here 
``superstandard'' means that we consecutively fill the first row,
followed by the second row, etc. 

Below, 
we observe there are only two tableaux $T,U$ 
on $\nu/\lambda$ that $K$-rectify to $S_{\mu}$:
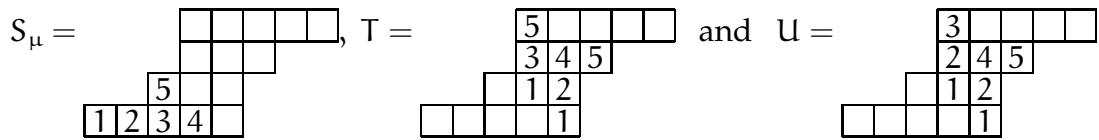
\[S_{\mu}=\tableau{&&&{\ }&{\ }&{\ }&{\  }&{\ }\\&&&{\ }&{\ }&{\ }\\&&{5 }&{\ }&{\ }\\
{1 }&{2 }&{3 }&{4 }&{\ }}, \ 
T=\tableau{&&&{5 }&{\ }&{\ }&{\  }&{\ }\\&&&{3 }&{4 }&{5}\\&&{\ }&{1 }&{2 }\\
{\ }&{\ }&{\ }&{\ }&{1 }} \mbox{ \ \ and \ \ }
U=\tableau{&&&{3 }&{\ }&{\ }&{\  }&{\ }\\&&&{2 }&{4 }&{5}\\&&{\ }&{1 }&{2 }\\
{\ }&{\ }&{\ }&{\ }&{1 }}
\]
Therefore, our conjecture states that: 
\[C_{(1,1,2,1),(1,1,2,1)}^{(1,1,2,4,3,1)}({\mathbb O}{\mathbb P}^2)=(-1)^{12-5-5}2=2.\]
The reader can check that the rectification order does not affect the
result. For either $T$ or $U$, there are three initial ways to begin the
${K}$-rectification, after which, all further $K{\tt jdt}$ slides
are forced. 
\end{Example}

Note that once one establishes an analogue of 
Theorem~\ref{thm:welldefined},
one can give an easy modification of the proof of 
associativity in Section~6 to 
establish that Conjecture~\ref{conjecture:mainone} 
defines an associative product. One can check that 
the analogue of Theorem~\ref{thm:welldefined} holds in specific instances, say,
with the help of a computer. Indeed, we have made exhaustive checks
when $G/P$ is the odd orthogonal Grassmannian 
$OG(5,11)$ and when it is the Cayley plane ${\mathbb O}{\mathbb P}^2$,
corresponding to the types $B_5$ and $E_6$. We
also made numerous checks in the case of the space $G_{\omega}({\mathbb O}^3, {\mathbb O}^6)$ associated to $E_7$, which while not exhaustive, left us 
convinced.

We emphasize that this rule agrees in type $A$
with the correct product, and as well as in cohomology
for all minuscule cases. We also have some computational evidence that our
numbers agree with small known cases of Schubert structure constants in type $B$
(as supplied to us by M.~Shimozono in private correspondence), although
admittedly this is not a convincing amount of evidence on its own. Part of
the difficulty in checking Conjecture~\ref{conjecture:mainone} 
is that it seems to be a challenging
task to construct efficient software to compute the $K$-theory 
Schubert structure constants for the main cases of the minuscule $G/P$'s
outside of type $A$. In principle, such an algorithm is 
linear algebra using torus-equivariant fixed-point 
localization methods such as \cite{Willems}.

Granted associativity, 
the conjectures would follow if they agree with multiplication in
$K(G/P)$ 
whenever $\mu$ is drawn from some
set of multiplicative generators ${\mathcal P}$ for $K(G/P)$. (That is,
they agree with a ``Pieri rule''.) We will report on our progress on 
these conjectures in forthcoming work.

We also mention that the results of Sections~2-4 also have straightforward
minuscule generalizations in cohomology, cf., \cite{Thomas.Yong:DE}.

\section{Counterexamples}
It is interesting that natural analogues of a 
number of results valid in the standard Young tableau theory
are actually false in our setting. We have
already seen in the introduction that in general $K{\tt rect}$ is not 
well-defined. This aspect can also be blamed for the following
two other situations where counterexamples exist:

\noindent
\emph{Haiman's dual equivalence:} 
One can define {\bf $K$-theoretic dual equivalence},
extending ideas in Haiman's \cite{Haiman:DE}. Two increasing tableaux
are {\bf $K$-dual equivalent} if any sequence of slides
$(\{x_{i_1}^{(1)}\},\ldots,\{x_{i_k}^{(k)}\})$ 
for $T$ and $U$ 
results in increasing tableaux of the
same shape. In this case we write $T\equiv_{D} U$. By definition, 
$T\equiv_{D} U$ implies ${\tt shape}(T)={\tt shape}(U)$. 

One application of this theory (in the classical setting) is that
it leads to a proof of the fundamental theorem of jeu de taquin. For
a minuscule (but not $K$-theoretic) generalization, see 
\cite{Thomas.Yong:DE}. However, it is important for this application
that all standard Young tableau of the same shape are dual equivalent.
In view of Theorem~\ref{thm:welldefined}, it is not surprising that this 
is not true in our setting. Consider the computations
\[K{\tt infusion}_2\left(\tableau{{1}&{3}\\{2}},\tableau{&&{2}\\&{1}&{4}\\{1}&{3}}\right)=\tableau{&&&\\&&{1}\\{2}&{3}} \mbox{\ \ and \ \ }
K{\tt infusion}_2\left(\tableau{{1}&{2}\\{3}},\tableau{&&{2}\\&{1}&{4}\\{1}&{3}}\right)=\tableau{&&&\\&{1}&{2}\\{1}&{3}}.
\]
These calculations represent two sequences of $K{\tt jdt}$ slides applied
to different tableaux of the same shape $(2,1)$, but whose results are
tableaux of different (skew) shapes.

\medskip
\noindent
\emph{Cartons:} In an earlier paper \cite{Thomas.Yong:IV}, we gave
an $S_3$-symmetric Littlewood-Richardson rule in terms of \emph{cartons}. This
idea also has a minuscule extension (which we will report on elsewhere). However,
the na\"{i}ve $K$-theoretic generalization does not work.

Briefly, the {\bf carton} of \cite{Thomas.Yong:IV} 
is a three-dimensional box with 
a grid drawn rectilinearly on the six faces of its surface, each of
whose sides are growth diagrams. We fix at
the outset standard Young tableaux of shape $\lambda,\mu$ and $\nu$
along three edges. Shapes are associated to each vertex so that the 
Fomin growth conditions (F1) and (F2) reproduced in Section~2 hold. The
number of such cartons (with fixed initial data) is equal to the classical
Littlewood-Richardson number. 

The temptation is to attempt to generalize this to $K$-theory by replacing the
initial standard Young tableau with superstandard tableau of shapes
$\lambda,\mu$ and $\nu$, and to instead utilize the growth conditions 
(G1) and (G2) we introduced in
Section~2. This does not work: 
one computes using Theorem~\ref{thm:KLR} that if $k=n-k=3$,
$\lambda=\mu=(2,1)$ and $\nu=(2)$ then
the constant $C_{(2,1),(2,1),(2)}:=C_{(2,1),(2,1)}^{(3,3,1)}=-2$.
However one cannot consistently complete a legal filling of this $K$-carton.

\begin{Remark}
These obstructions are closely related to failure of associativity 
of a certain tableau product defined in \cite[Section~3.7]{BKSTY}.
\end{Remark}

\section{Concluding remarks}

\subsection{Proctor's $d$-complete posets}
Proctor~\cite{Proctor} has studied the class of \emph{$d$-complete
posets}. These posets generalize those required in our
discussion of minuscule $G/P$ Schubert calculus; see 
also~\cite{Thomas.Yong, Thomas.Yong:DE}. In particular, $d$-complete
posets were shown by Proctor to have a well-defined {\it jeu de
taquin} procedure.  

It would be interesting to generalize our arguments 
to show that for any $d$-complete poset $D$, there is an associative
ring $K(D)$ with an additive ${\mathbb Z}$-basis indexed by lower
order ideals of $D$ and structure constants defined by a 
rule generalizing Theorem~\ref{thm:KLR}. Observing that our notions of
$K{\tt jdt}, K{\tt rect}$ \emph{a priori} make sense in this more general context, we ask:

\begin{Problem}
\label{conj:proctor}
Fix a $d$-complete poset. For which classes of tableaux
${\mathcal C}=\{C_{\mu}\}$ (indexed by lower order ideals $\mu$ of $D$) 
is it true that an analogue Theorem~\ref{thm:welldefined} holds (that is
if $K{\tt rect}(T)=C\in {\mathcal C}$ under one rectification order,
this holds for any rectification order)?
\end{Problem} 

It seems plausible that good classes ${\mathcal C}$ that play the role
of the superstandard tableaux of Theorem~\ref{thm:welldefined} always exist.
As we have said, for the minuscule cases, we believe that the
superstandard tableaux suffice. Perhaps this also holds more generally.

Assuming the Conjecture holds, 
one would also like to find a geometric origin to the ring $K(D)$ 
(outside of the cases where it should be isomorphic to the $K$-theory
ring of a minuscule $G/P$).

\subsection{A product-differences conjecture}

Let $\lambda,\mu\in {\mathbb Y}$. Since this poset is in fact a lattice,
we can speak of their {\bf meet} $\lambda\wedge\mu$ and {\bf join}
$\lambda\vee\mu$. 

\begin{Conjecture}
\label{conj:proddiff}
Suppose $\lambda,\mu \subseteq \Lambda$.
Let 
\[[{\mathcal O}_{X_{\lambda\wedge\mu}}][{\mathcal O}_{X_{\lambda\vee\mu}}]-
[{\mathcal O}_{X_{\lambda}}][{\mathcal O}_{X_{\mu}}]=\sum_{\nu}d_{\nu}[{\mathcal O}_{X_{\nu}}].\]
Then $(-1)^{|\nu|-|\lambda|-|\mu|}d_{\nu}\geq 0$.
\end{Conjecture}

This conjecture generalizes a theorem in the cohomological case due
to \cite{bindiff}; see related work by \cite{Okounkov, FFLP, Chindris}.
(We also know of no counterexample for the corresponding 
minuscule conjecture, even in the cohomology case.)

\begin{Example}
Let $\lambda=(4,2,1),\mu=(3,3,2)\subseteq \Lambda=4\times 5$. The join
is the unique minimal shape that contains $\lambda$ and $\mu$, i.e., 
$\lambda\vee\mu=(4,3,2)$. Similarly, the meet is the unique maximal
shape contained in $\lambda$ and $\mu$. Hence $\lambda\wedge\mu=(3,2,1)$. One
computes using Theorem~\ref{thm:KLR} (or otherwise), preferably with the
help of computer, that:
\begin{multline}
[{\mathcal O}_{X_{(4,3,2)}}]\cdot [{\mathcal O}_{X_{(3,2,1)}}]
-[{\mathcal O}_{X_{(4,2,1)}}]\cdot [{\mathcal O}_{X_{(3,3,2)}}]
=([{\mathcal O}_{X_{(5,5,3,2)}}]+2[{\mathcal O}_{X_{(5,5,4,1)}}]
+[{\mathcal O}_{X_{5,5,5}}]+[{\mathcal O}_{X_{5,4,4,2}}])\\ \nonumber
-(3[{\mathcal O}_{X_{(5,5,5,1)}}]+[{\mathcal O}_{X_{(5,5,3,3)}}]
+5[{\mathcal O}_{X_{(5,5,4,2)}}]+[{\mathcal O}_{X_{(5,4,4,3)}}])
+(3[{\mathcal O}_{X_{(5,5,5,2)}}]+3[{\mathcal O}_{X_{(5,5,4,3)}}])\\ \nonumber
-([{\mathcal O}_{X_{(5,5,5,3)}}]),
\end{multline}
in agreement with Conjecture~\ref{conj:proddiff}.
\end{Example}

%
%

\subsection{Hecke insertion and factor sequence formulae}

In \cite{BKSTY} a generalizaion of the Robinson-Schensted and
Edelman-Greene insertion algorithms was given. In fact, increasing tableaux
also play a prominent role there, although in a different, but 
related way. As we have mentioned in the introduction, this
will be explored, in part, in \cite{Thomas.Yong:VI}, in connection to 
longest strictly increasing subsequences in random words. There we show 
that the insertion tableau of a word under Hecke insertion can be 
alternatively computed as a $K$-rectification of a permutation
tableau (for a particular
choice of rectification order).

 Another sample question: is there a ``plactification map'' in
the sense of \cite{Reiner.Shimozono}? 

We believe that further developing this connection may
allow one to, for example, prove a $K$-theory analogue of the 
``factor sequence formula'' conjectured in \cite{Buch.Fulton} and 
proved in \cite{KMS}, which is a problem that has remained open in this
topic, see \cite{Buch:original, Buch:JAMS}. (In \cite{BKSTY} a different factor sequence
formula, generalizing the one given in \cite{Buch:JAMS}, was given.)

\section*{Appendix: Grothendieck polynomials}

The goal of this appendix is to provide combinatorial background
for the results of Sections 1--7, in terms of the Grothendieck polynomials of
Lascoux and Sch\"{u}tzenberger \cite{LS:Hopf}. This presentation is
not needed for the paper.

Fix a shape $\lambda$ and define a {\bf set-valued tableau} $T$ 
to be an assignment of nonempty sets of natural numbers to each box of
$\lambda$ \cite{Buch:KLR}. Such a tableau is {\bf semistandard} if for
every box, the largest entry is weakly smaller than the minimum entry of
the box immediately to its right and strictly smaller than the minimum
entry of the box immediately below it. The {\bf ordinary} case is when $T$
assigns a singleton to each box.
\begin{figure}[h]
\[T_1=\tableau{{1}&{2}&{4}&{4}&{6}\\{2}&{3}&{5}\\{4}}
\ \ T_2=\ktableau{{1,2}&{2,3}&{4,5,6}&{6,7}&{7,8}\\{3,4}&{4,5}&{7}\\{6,7,8}}\]
\caption{An ordinary and a set-valued semistandard tableau}
\end{figure}

Associate to each semistandard tableau a {\bf weight} 
\[\omega(T):=(-1)^{|T|-|\lambda|}{\bf x}^{T}\] 
where here ${\bf x}^{T}=x_1^{i_1}x_2^{i_2}\cdots$ if
$i_j$ is the number of $j$'s appearing in $T$, and $|T|$ is the number of
entries of $T$. For example, we have 
\[\omega(T_1)=x_1 x_2^2 x_3 x_4^3 x_5 x_6 \mbox{ \ and \ }
\omega(T_2)=(-1)^{19-9}x_1 x_2^2 x_3^2 x_4^3 x_5^2 x_6^3 x_7^4 x_8^2.\]

The {\bf Grothendieck polynomial} is defined as
\[G_{\lambda}(x_1, x_2,\ldots, x_k):=\sum_{T} {\omega}(T)\]
with the sum over all set-valued semistandard tableaux using the labels
of size at most $k$. This is an
inhomogenous symmetric polynomial whose lowest degree ($=|\lambda|$)
homogeneous component is equal to the {\bf Schur polynomial}
$s_{\lambda}(x_1,x_2,\ldots x_k)$.

It is not immediately obvious from the definitions, but true
\cite{Buch:KLR} (for an alternative proof, see \cite{BKSTY}) that the $G_{\lambda}(x_1,\ldots,x_k)$ (for $\lambda$ with at most $k$ parts) form a ${\mathbb Z}$-linear basis for the ring of symmetric polynomials
in $x_1,\ldots, x_k$ (say, with coefficients in ${\mathbb Q}$). Thus we can
write
\[G_{\lambda}(x_1,\ldots,x_k)G_{\mu}(x_1,\ldots,x_k)=\sum_{\nu}C_{\lambda,\mu}^{\nu}G_{\nu}(x_1,\ldots,x_k).\]
The coefficients $C_{\lambda,\mu}^{\nu}$ agree with the $K$-theory
structure constants for $Gr(k,{\mathbb C}^n)$ whenever $\nu\subseteq \Lambda$.

There are more general Grothendieck polynomials ${\mathfrak G}_{\pi}(x_1,\ldots,x_n)$ defined
in \cite{LS:Hopf} for any permutation $\pi\in S_n$. The polynomials
$G_{\lambda}$ amount to the case that $\pi$ is {\bf Grassmannian}: it has
a unique descent at position $k$. In \cite{BKTY:groth} a formula was first 
given that expresses any ${\mathfrak G}_{\pi}$ in terms of the $G_{\lambda}$'s.
Other formulas for both 
${\mathfrak G}_{\pi}$ and $G_{\lambda}$ are also available, see, 
e.g.,~\cite{BKSTY, KY, KMY, lascoux:transition} and the references therein.

\section*{Acknowledgments}
This work was partially 
completed while HT was visiting the Norges \linebreak Teknisk-Naturvitenskapelige
Universitet; he would like to thank the 
Institutt for Matematiske Fag for its hospitality.  
AY utilized the resources of
the Fields Institute, Toronto, while a visitor there. 
We thank Allen Knutson, Victor Reiner and Mark Shimozono for helpful 
discussions, as well as Anders Buch for supplying us with software to 
independently compute the $K$-theoretic numbers $C_{\lambda,\mu}^{\nu}$
(for Grassmannians).

\end{document}